# RANK-3 STABLE BUNDLES ON RATIONAL RULED SURFACES

WEI-PING LI AND ZHENBO QIN

ABSTRACT. In this paper, we compare the moduli spaces of rank-3 vector bundles stable with respect to different ample divisors over rational ruled surfaces. We also discuss the irreducibility, unirationality, and rationality of these moduli spaces.

## INTRODUCTION.

This paper is the continuation of our previous paper [9] on higher rank stable bundles over algebraic surfaces. One of the main results in [9] is that on a rational ruled surface $X$, the moduli space of vector bundles stable with respect to some *suitable* ample divisor is irreducible and unirational whenever the moduli space is nonempty. It is also announced in [9] that every irreducible component of an arbitrary moduli space is unirational. Therefore, it is very natural to ask whether an arbitrary moduli space is irreducible. In this paper, we shall answer this question for rank-3 stable bundles on a rational ruled surface $X$.

Our idea is the following. We start with a *suitable* ample divisor $L_0$ on $X$ such that the moduli space of $L_0$-stable bundles is irreducible and unirational whenever it is nonempty. For an arbitrary ample divisor $L$ on $X$, we compare the two moduli spaces of rank-3 bundles stable with respect to $L_0$ and $L$ respectively. We generalize the methods in [11, 12]. It follows from the Harder-Narasimhan filtration that every rank-3 bundle which is contained in the difference of the two moduli spaces can be written either as an extension of a rank-2 torsion free sheaf by a line bundle, or as an extension of a rank-1 torsion free sheaf by a rank-2 vector bundle; in either case, the extension relates to a wall in the Kähler cone $\mathbb{C}_X$ of $X$, which is the closed convex cone in $\mathrm{Num}(X) \otimes \mathbb{R}$ spanned by (the images of) ample divisors. Using some techniques similar to those in [11], we estimate the number of moduli of the rank-3 vector bundles in the difference of the two moduli spaces. ¿From these estimates, we are able to draw conclusions about the irreducibility and unirationality of the moduli space of rank-3 $L$-stable bundles on $X$.

Next, we state our main results in this paper. Let $X$ be a rational ruled surface, and let $\pi : X \to \mathbb{P}^1$ be the ruling. Let $f$ be a fiber of $\pi$, and let $\sigma$ be a section to $\pi$ such that $\sigma^2$ is the least. Put $e = -\sigma^2$. Let $c_1$ be a divisor, $c_2$ be an integer, and $L$ be an ample divisor on $X$. We use $\mathfrak{M}_L(3; c_1, c_2)$ to stand for the moduli space of rank-3 $L$-stable vector bundles on $X$ with the first Chern class $c_1$ and the second Chern class $c_2$. It is well-known (see [8]) that $\sigma$ and $f$ generates the Picard group of $X$. If $L = (x\sigma + yf)$, then we let $r_L = y/x$. We have the following two results for $c_1 = 0$ and $c_1 = (\sigma + f)$ respectively.

**Theorem A.** (i) The moduli space $\mathfrak{M}_L(3; 0, c_2)$ is nonempty if and only if $c_2 \geq 3$;

(ii) If $c_2 \geq 3$, then $\mathfrak{M}_L(3; 0, c_2)$ is smooth, irreducible and unirational with dimension $(6c_2 - 8)$.

**Theorem B.** Assume that the integer $c_2$ is even.

(i) If $c_2 < 2$, then the moduli spaces $\mathfrak{M}_L(3; \sigma + f, c_2)$ is empty;

(ii) Assume that $c_2 \geq 2$. If $r_L \geq e + (3c_2 - 2)/2$ or $r_L \leq 2/(3c_2 - 2)$, then $\mathfrak{M}_L(3; \sigma + f, c_2)$ is empty. If $2/(3c_2 - 2) < r_L < e + (3c_2 - 2)/2$, then $\mathfrak{M}_L(3; \sigma + f, c_2)$ is irreducible, smooth, and unirational with dimension $(6c_2 + 2e - 12)$.

Assume that the integer $c_2$ is odd.

(i) If $c_2 < 3$, then the moduli space $\mathfrak{M}_L(3; \sigma + f, c_2)$ is empty;

(ii) Assume that $c_2 \geq 3$. If $r_L \geq e + (3c_2 - 5)/2$ or $r_L \leq 2/(3c_2 - 5)$, then $\mathfrak{M}_L(3; \sigma + f, c_2)$ is empty. If $2/(3c_2 - 5) < r_L < e + (3c_2 - 5)/2$, then $\mathfrak{M}_L(3; \sigma + f, c_2)$ is irreducible, smooth, and rational with dimension $(6c_2 + 2e - 12)$.

We notice that the blow-up of the projective plane $\mathbb{P}^2$ at one point is the rational ruled surface with $e = 1$. Let $L$ be the divisor represented by a line in $\mathbb{P}^2$. In [3], Drezet and Le Potier have determined the nonemptyness of the moduli space $\mathfrak{M}_L(3; c_1, c_2)$, and showed that $\mathfrak{M}_L(3; c_1, c_2)$ is irreducible whenever it is nonempty. In particular, by the Theorem B in [3], $\mathfrak{M}_L(3; 0, c_2)$ is nonempty if and only if $c_2 \geq 3$, and $\mathfrak{M}_L(3; L, c_2)$ is nonempty if and only if $c_2 \geq 2$. As an application of Theorem A and Theorem B, we improve upon Drezet and Le Potier's result when $r = 3$.

**Theorem C.** If $c_2 \geq 3$, then $\mathfrak{M}_L(3; 0, c_2)$ is irreducible, smooth, and unirational with dimension $(6c_2 - 8)$. If $c_2 \geq 2$, then $\mathfrak{M}_L(3; L, c_2)$ is irreducible, smooth, and unirational (rational when $c_2$ is odd) with dimension $(6c_2 - 10)$.

We remark that while we consider the case of rank-3 stable vector bundles only, we don't loss the essence of higher rank stable bundles and at the same time calculations are much easier. The reason that we choose our surfaces to be rational ruled surfaces is that the Neron-Severi group of these surfaces is very simple. We also mention that while we are writing up this paper, we noticed a recent paper by Walter [13], which covers most part of our results. However, our results are more explicit, and our approach should be applicable to other algebraic surfaces such as algebraic surfaces with Kodaira dimension zero.

Finally, as a convention, we make no distinction between a locally free sheaf and the corresponding vector bundle. Our paper is organized as follows. In section 1, we compare stabilities with respect to different ample divisors. In section 2, we estimate the number of moduli of rank-3 bundles in the difference of two moduli spaces. In section 3, we prove the three theorems above.

## 1. Comparison of stabilities.

In this section, $X$ stands for an arbitrary algebraic surface. We shall use the Harder-Narisimhan filtration to compare the stabilities with respect to different ample divisors. Then, we introduce walls and chambers, and study the stability of certain vector bundles which are given by some extensions related to walls.

### 1.1. Comparison of stabilities.

In this subsection, we recall definition of the stability in the sense of Mumford and Takemoto, and the Harder-Narasimhan filtration of a vector bundle with respect to an ample divisor. Then, we compare the stabilities with respect to different ample divisors on the surface $X$.

**Definition 1.1.** Let $L$ be an ample divisor on $X$, and let $V$ be a torsion free coherent sheaf. Put

$$\mu_L(V) = \frac{c_1(V) \cdot L}{\mathrm{rank}(V)}.$$

Then, $V$ is defined to be *L-stable* (respectively, *L-semistable*) if for any proper lower rank subsheaf $F$ of $V$, we have

$$\mu_L(F) < \mu_L(V) \ (\text{respectively, } \mu_L(F) \leq \mu_L(V)).$$

moreover, $V$ is defined to be *strictly L-semistable* if it is $L$-semistable but not $L$-stable, and to be *L-unstable* if it is not $L$-semistable. We will use $\mathfrak{M}_L(r; c_1, c_2)$ to represent the moduli space of rank-$r$ $L$-stable vector bundles $V$ with the first Chern class $c_1$ and the second Chern class $c_2$.

**Proposition 1.2.** *(see [7, 4]) Fix an ample divisor $L$ on $X$. Then, every bundle $V$ has a canonical (Harder-Narasimhan) filtration*

$$0 \subset V_1 \subset V_2 \subset \ldots \subset V_k \subset V$$

*such that $V_i$ is a proper subbundle of $V$ for all $i$, $V/V_i$ is torsion free, and $V_i/V_{i-1}$ is an $L$-semistable subsheaf of $V/V_{i-1}$ with the following two properties:*

(i) $\mu_L(V_i/V_{i-1}) = N_i = \max\{\mu_L(U)|\ U$ *is a subsheaf of $V/V_{i-1}\}$;*
(ii) $\mathrm{rank}(V_i/V_{i-1}) = \max\{\mathrm{rank}(U)|\ U$ *is a subsheaf of $V/V_{i-1}$ and $\mu_L(U) = N_i\}$.*

**Lemma 1.3.** *Let $L_1$ and $L_2$ be two ample divisors on $X$. Assume that $V$ is a $L_1$-stable rank-$r$ bundle and that $V$ is $L_2$-unstable with the Harder-Narasimhan filtration: $0 \subset V_1 \subset V$. Then, there exists an integer $i$ with $0 < i < r$ and a divisor $F$ such that*

(i) $(rF - ic_1(V)) \cdot L_1 < 0 < (rF - ic_1(V)) \cdot L_2$;
(ii) $-i(r-i)(r-1) \cdot [2rc_2(V)/(r-1) - c_1(V)^2] \leq (rF - ic_1(V))^2 < 0$.

*Proof.* Let $i$ and $F$ be the rank and the first Chern class of $V_1$ respectively.

(i) Since $\mu_{L_1}(V_1) < \mu_{L_1}(V)$ and $\mu_{L_2}(V_1) > \mu_{L_2}(V)$, the inequalities follow.

(ii) By (i), $(rF - ic_1(V)) \cdot L = 0$ for some ample divisor $L$. By the Hodge Index Theorem (see [8]), $(rF - ic_1(V))^2 < 0$. Next, we notice that both $V_1$ and $V/V_1$ are $L_2$-semistable. By the Bogomolov's inequality (see [1] ), we have

$$c_2(V_1) \geq \frac{i-1}{2i} \cdot F^2 \text{ and } c_2(V/V_1) \geq \frac{(r-i)-1}{2(r-i)} \cdot (c_1(V) - F)^2.$$

Since $c_2(V) = c_2(V_1) + c_2(V/V_1) + F \cdot (c_1(V) - F)$, we have

$$\begin{aligned} c_2(V) &\geq \frac{i-1}{2i} \cdot F^2 + \frac{(r-i)-1}{2(r-i)} \cdot (c_1(V) - F)^2 + F \cdot (c_1(V) - F) \\ &= \frac{r-1}{2r} \cdot c_1(V)^2 - \frac{(rF - ic_1(V))^2}{2i(r-i)r}. \end{aligned}$$

Rewriting this, we obtain the desired inequality. $\square$

**Lemma 1.4.** *Let $L_1$ and $L_2$ be two ample divisors on $X$. Assume that $V$ is a $L_1$-stable rank-$r$ bundle and that $V$ is $L_2$-unstable with the Harder-Narasimhan filtration: $0 \subset V_1 \subset V_2 \subset V$. Then, there exists an integer $i$ with $0 < i < r$ and a divisor $F$ such that*

(i) $(rF - ic_1(V)) \cdot L_1 < 0 < (rF - ic_1(V)) \cdot L_2$;
(ii) $-i(r-i)(r-1) \cdot [2rc_2(V)/(r-1) - c_1(V)^2] < (rF - ic_1(V))^2 < 0$.

*Proof.* Let $c_1 = c_1(V)$, $F_j = c_1(V_j)$ and $r_j = \text{rank}(V_j)$. Since $V$ is $L_1$-stable,

$$(rF_j - r_jc_1) \cdot L_1 < 0.$$

Choose a number $x_2$ with $0 < x_2 < (r_2 - r_1)/r_2$. Put

$$x_1 = \frac{(r_2 - r_1) - r_2x_2}{r_1} \text{ and } x = \frac{(r - r_2) + rx_1}{r_2}.$$

Then, $x_1 > 0$ and $x > 0$. Moreover, one checks the following equality:

$$[(r_2 - r)F_1 + (r - r_1)F_2 + (r_1 - r_2)c_1] + x(r_2F_1 - r_1F_2) = x_1(rF_1 - r_1c_1) + x_2(rF_2 - r_2c_1).$$

Thus, either $[(r_2 - r)F_1 + (r - r_1)F_2 + (r_1 - r_2)c_1] \cdot L_1 < 0$ or $(r_2F_1 - r_1F_2) \cdot L_1 < 0$. In the following, we consider these two cases separately.

Case 1: $[(r_2 - r)F_1 + (r - r_1)F_2 + (r_1 - r_2)c_1] \cdot L_1 < 0$. In this case, let $i = r_1$ and $F = F_1$. Since $V$ is $L_1$-stable, $(rF - ic_1) \cdot L_1 < 0$. Since $0 \subset V_1 \subset V_2 \subset V$ is the Harder-Narasimhan filtration of $V$ with respect to $L_2$, $\mu_{L_2}(V_1) > \mu_{L_2}(V)$, so $(rF - ic_1) \cdot L_2 > 0$. Thus, (i) and the second inequality in (ii) follow by the similar argument in the proof of Lemma 1.3. Next, we prove the first inequality in (ii). Consider the two exact sequences:

$$0 \to V_1 \to V \to V/V_1 \to 0 \qquad \text{and} \qquad 0 \to V/V_1 \to V/V_2 \to 0. \qquad (1.5)$$

Put $G = [(r_2 - r)F_1 + (r - r_1)F_2 + (r_1 - r_2)c_1]$. By Proposition 1.2,

$$\mu_{L_2}(V_2/V_1) > \mu_{L_2}(V/V_1).$$

This implies that $G \cdot L_2 > 0$. Since $G \cdot L_1 < 0$, $G^2 < 0$ by the Hodge index theorem. Note that $V_1, V_2/V_1$ and $V/V_2$ are all $L_2$-semistable. Applying the method as in the proof of Lemma 1.3 (ii) to the second exact sequence in (1.5), we conclude that

$$c_2(V/V_1) \geq \frac{(r-r_1)-1}{2(r-r_1)}(c_1 - F_1)^2 - \frac{G^2}{2(r_2 - r_1)(r - r_2)(r - r_1)}$$

$$> \frac{(r-r_1)-1}{2(r-r_1)}(c_1 - F_1)^2.$$

Since $V_1$ is $L_1$-semistable, again by Bogolomov's inequality, we have

$$c_2(V_1) \geq \frac{r_1 - 1}{2r_1} F^2.$$

By the similar argument as in the proof of Lemma 1.3 (ii), we get

$$c_2(V) > \frac{r-1}{2r} \cdot c_1^2 - \frac{(rF_1 - r_1 c_1)^2}{2rr_1(r - r_1)}.$$

Therefore, $-i(r-i)(r-1) \cdot [2rc_2(V)/(r-1) - c_1(V)^2] < (rF - ic_1(V))^2$.

Case 2: $(r_2 F_1 - r_1 F_2) \cdot L_1 < 0$. In this case, let $i = r_2$ and $F = F_2$. Consider the two exact sequences:

$$0 \to V_2 \to V \to V/V_2 \to 0 \qquad \text{and} \qquad 0 \to V_1 \to V_2 \to V_2/V_1 \to 0.$$

A similar argument as in Case 1 shows the desired inequality. $\square$

**Theorem 1.6.** *Let $L_1$ and $L_2$ be two ample divisors on $X$. Assume that $V$ is a rank-3 bundle $L_1$-stable but $L_2$-unstable. Then, there exists an integer $i$ with $0 < i < 3$ and a divisor $F$ such that*

  (i) $(3F - ic_1(V)) \cdot L_1 < 0 < (3F - ic_1(V)) \cdot L_2$;
  (ii) $-4 \cdot [3c_2(V) - c_1(V)^2] \leq (3F - ic_1(V))^2 < 0$.

*Proof.* Assume that the Harder-Narasimhan filtration of $V$ with respect to $L_2$ is:

$$0 \subset V_1 \subset \ldots \subset V_k \subset V$$

where $k = 1$ or $2$. If $k = 1$, then the conclusions follow from Lemma 1.3; if $k = 2$, then the conclusions follow from Lemma 1.4. $\square$

*Remark* 1.7. Assume that $V$ is a rank-3 bundle $L_1$-stable but $L_2$-unstable and that the Harder-Narasimhan filtration of $V$ with respect to $L_2$ is:

$$0 \subset V_1 \subset \ldots \subset V_k \subset V$$

¿From the above proofs, we see that one of the following is true:

  (a) $k = 1$, $F = c_1(V_1)$, and $i = \operatorname{rank}(V_1) = 1$;
  (b) $k = 1$, $F = c_1(V_1)$, and $i = \operatorname{rank}(V_1) = 2$;
  (c) $k = 2$, $F = c_1(V_2)$, and $i = \operatorname{rank}(V_2) = 2$;
  (d) $k = 2$, $F = c_1(V_1)$, and $i = \operatorname{rank}(V_1) = 1$.

We shall analysis these four cases in the next section.

## 1.2. Walls and chambers.

Let $\mathbb{C}_X$ be the Kähler cone of $X$, which is the closed convex cone in $\mathrm{Num}(X) \otimes \mathbb{R}$ spanned by ample divisors. In this subsection, we define walls and chambers in $\mathbb{C}_X$.

Fix a divisor $c_1$ on $X$ and an integer $c_2$ such that

$$\frac{2r}{r-1} \cdot c_2 - c_1^2 > 0.$$

The following generalizes the Definition 1.1 in [11].

**Definition 1.8.** (i) For $\zeta \in \mathrm{Num}(X) \times \mathbb{R}$, we define $W^\zeta$ to be the set

$$W^\zeta = \mathbb{C}_X \cap \{ h \in \mathrm{Num}(X) \times \mathbb{R} | \zeta \cdot h = 0 \};$$

(ii) Fix a positive integer $r$. We define $\mathcal{W}_{r;c_1,c_2}$ to be the set whose elements are $W^\zeta$, where $\zeta$ is the numerical equivalence class of $(rF - ic_1)$ for some divisor $F$ and $0 < i < r$, and satisfies

$$-i(r-i)(r-1) \cdot [\frac{2r}{r-1} \cdot c_2 - c_1^2] \leq \zeta^2 < 0;$$

(iii) *A wall of type* $(r; c_1, c_2)$ is an element in $\mathcal{W}_{r;c_1,c_2}$. *A chamber of type* $(r; c_1, c_2)$ is a connected component in the complement $\mathbb{C}_X - \mathcal{W}_{r;c_1,c_2}$.

By the results in Chapter II of [5] or the Proposition 2.1.6 in Chapter I of [12], the set $\mathcal{W}_{r;c_1,c_2}$ of walls of a fixed type is locally finite. In view of this, we make the following definition.

**Definition 1.9.** Let $\mathcal{C}$ be a chamber of type $(r; c_1, c_2)$, and let $\overline{\mathcal{C}}$ be the closure of $\mathcal{C}$ in $\mathbb{C}_X$. Then, *a face of* $\mathcal{C}$ or $\overline{\mathcal{C}}$ is a codimension-1 subset of $\overline{\mathcal{C}}$ of the form $W \cap \overline{\mathcal{C}}$ where $W$ is a wall of type $(r; c_1, c_2)$.

*Remark 1.10.* (i) If $\mathcal{F}$ is a face of $\mathcal{C}$, then there is unique chamber $\mathcal{C}' \neq \mathcal{C}$ which also has $F$ as a face. In this case, the chamber $\mathcal{C}$ and $\mathcal{C}'$ lie on opposite sides of the wall $W$ containing the common face $F$.

(ii) We explain the geometric meaning of Theorem 1.6. Let notations be as in Theorem 1.6. Put $\zeta = (3F - ic_1(V))$, and $L = -(\zeta \cdot L_1)L_2 + (\zeta \cdot L_2)L_1$. Then, $\zeta$ defines a wall $W^\zeta$ of type $(3; c_1, c_2)$, and $L \in W^\zeta$. Moreover, the wall $W^\zeta$ separates the two ample divisors $L_1$ and $L_2$.

**Corollary 1.11.** *Let* $L_1$ *and* $L_2$ *be two ample divisors in the same chamber* $\mathcal{C}$ *of type* $(3; c_1, c_2)$. *Then,* $\mathfrak{M}_{L_1}(3; c_1, c_2)$ *and* $\mathfrak{M}_{L_2}(3; c_1, c_2)$ *are the same.*

*Proof.* Follows immediately from Theorem 1.6 and Remark 1.10 (ii). $\square$

In the next lemma, we study the relation between the strict semistability and ample divisors which do not lie on any wall.

**Lemma 1.12.** *Let* $L$ *be an ample divisor which does not lie on any wall of type* $(r; c_1, c_2)$. *Let* $V$ *be a strictly* $L$-*semistable vector bundle over* $X$ *with* $c_1(V) = c_1$ *and* $c_2(V) = c_2$. *Then,* $V$ *sits in an exact sequence*

$$0 \to V_1 \to V \to V_2 \to 0 \tag{1.13}$$

where $0 < \text{rank}(V_1) < r$ and $[r \cdot c_1(V_1) - \text{rank}(V_1) \cdot c_1] \equiv 0$. In particular, when $r = 2$ or $3$, we have $c_1 \equiv 0 \pmod{r}$.

*Proof.* Since $V$ is strictly $L$-semistable, there exists a nonzero proper locally free subsheaf $V_1$ of $V$ such that $V/V_1$ is torsion free and that $\mu_L(V_1) = \mu_L(V)$. Putting $V_2 = V/V_1$, we obtain the exact sequence (1.13). For simplicity, let $c_1' = c_1(V_1)$ and $r' = \text{rank}(V_1)$. Since $\mu_L(V_1) = \mu_L(V)$, we have $(r \cdot c_1' - r' \cdot c_1) \cdot L = 0$. Thus, by the Hodge index theorem, $(r \cdot c_1' - r' \cdot c_1)^2 \leq 0$ with equality if and only if $(r \cdot c_1' - r' \cdot c_1) \equiv 0$. Thus, it suffices to show that $(r \cdot c_1' - r' \cdot c_1)^2$ can not be negative. Assume that $(r \cdot c_1' - r' \cdot c_1)^2 < 0$. We want to draw a contradiction. First of all, we notice that $\mu_L(V_1) = \mu_L(V) = \mu_L(V_2)$; since $V$ is $L$-semistable, both $V_1$ and $V_2$ must be $L$-semistable. By the same argument as in the proof of Lemma 1.3 (ii), we conclude that

$$-r'(r - r')(r - 1)\left(\frac{2r}{r-1} \cdot c_2 - c_1^2\right) \leq (r \cdot c_1' - r' \cdot c_1)^2.$$

Thus, $(r \cdot c_1' - r' \cdot c_1)$ defines a wall of type $(r; c_1, c_2)$; moreover, since $(r \cdot c_1' - r' \cdot c_1) \cdot L = 0$, the wall $W^{(r \cdot c_1' - r' \cdot c_1)}$ contains $L$. But this contradicts to our assumption that $L$ does not lie on any wall of type $(r; c_1, c_2)$. $\quad\square$

**Proposition 1.14.** *Let $L_1$ and $L_2$ be two ample divisors. Assume that $V \in \mathfrak{M}_{L_1}(r; c_1, c_2)$ and that $L_2$ does not lie on any wall of type $(r; c_1, c_2)$. Then, $V$ can not be strictly $L_2$-semistable.*

*Proof.* Otherwise, by Lemma 1.12, $V$ sits in an exact sequence

$$0 \to V_1 \to V \to V_2 \to 0$$

where $0 < \text{rank}(V_1) < r$ and $[r \cdot c_1(V_1) - \text{rank}(V_1) \cdot c_1] \equiv 0$. But then $\mu_{L_1}(V_1) = \mu_{L_1}(V)$, contradicting to the assumption that $V$ is $L_1$-stable. $\quad\square$

**Corollary 1.15.** *Let $L_1$ and $L_2$ be two ample divisors. Assume that $L_1$ lies in a chamber $\mathcal{C}$ of type $(3; c_1, c_2)$ and that $L_2$ lies in the closure $\overline{\mathcal{C}}$ of $\mathcal{C}$. Then, $\mathfrak{M}_{L_2}(3; c_1, c_2)$ is naturally embedded in $\mathfrak{M}_{L_1}(3; c_1, c_2)$.*

*Proof.* When $L_2 \in \mathcal{C}$, then $\mathfrak{M}_{L_2}(3; c_1, c_2)$ is naturally identified with $\mathfrak{M}_{L_1}(3; c_1, c_2)$ by Theorem 1.6. Assume that $L_2 \in (\overline{\mathcal{C}} - \mathcal{C})$. Let $V \in \mathfrak{M}_{L_2}(3; c_1, c_2)$. By Theorem 1.6, $V$ must be $L_1$-semistable: otherwise, $L_1$ and $L_2$ would be separated by a wall. Since $L_1$ lies in a chamber and $V$ is $L_2$-stable, by Proposition 1.14, $V$ can not be strictly $L_1$-semistable; thus, $V$ is $L_1$-stable. It follows that $\mathfrak{M}_{L_2}(3; c_1, c_2)$ is naturally embedded in $\mathfrak{M}_{L_1}(3; c_1, c_2)$. $\quad\square$

## 1.3. Stable rank-3 bundles given by certain extensions.

To start with, we recall a result in [12]. Let $\zeta = (2F - c_1)$ define a nonempty wall of type $(2; c_1, c_2)$, and let $\mathcal{C}$ be a chamber of type $(2; c_1, c_2)$ such that its closure $\overline{\mathcal{C}}$ intersects with $W^\zeta$ and that $L_1 \cdot \zeta < 0$ for $L_1 \in \mathcal{C}$. Assume that $V$ with $c_1(V) = c_1$ and $c_2(V) = c_2$ is a rank-2 bundle given by a nontrivial extension

$$0 \to \mathcal{O}_X(F) \to V \to V' \to 0.$$

Then, by the Theorem 1.2.3 in Chapter II of [12], $V$ must be $L_1$-stable. The goal of this subsection is to prove a similar result for the case of rank-3. Our result says that rank-3 bundles given by certain extensions which are related to walls must be stable with respect to some ample divisor.

**Proposition 1.16.** *Let* $\zeta = (3F - c_1)$ *define a nonempty wall of type* $(3; c_1, c_2)$, *and let* $\mathcal{C}$ *be a chamber of type* $(3; c_1, c_2)$ *such that its closure* $\overline{\mathcal{C}}$ *intersects with* $W^{\zeta}$ *and that* $L_1 \cdot \zeta < 0$ *for* $L_1 \in \mathcal{C}$. *Let* $L_0$ *be an ample divisor contained in* $(\overline{\mathcal{C}} \cap W^{\zeta})$. *Assume that* $V$ *with* $c_1(V) = c_1$ *and* $c_2(V) = c_2$ *is a rank-3 bundle given by a nontrivial extension*

$$0 \to \mathcal{O}_X(F) \to V \to V' \to 0 \tag{1.17}$$

*where* $V'$ *is* $L_0$-*semistable such that no* $L_0$-*destablizing rank-1 subsheaf can be lifted to a subsheaf of* $V$. *If* $c_1 \not\equiv 0 \,(\mathrm{mod}\,3)$ *and if* $V'$ *is* $L_1$-*semistable, then* $V$ *is* $L_1$-*stable.*

*Proof.* Since $c_1 \not\equiv 0 \,(\mathrm{mod}\,3)$ and since $L_1$ does not lie on any wall of type $(3; c_1, c_2)$, by Lemma 1.12, it suffices to show that $V$ is $L_1$-semistable. Assume the contrary, that is, $V$ is $L_1$-unstable. We shall draw a contradiction. First of all, we prove the following claim.

**Claim.** *Assume that* $U$ *is a locally free subsheaf of* $V$ *with* $0 < \mathrm{rank}(U) < 3$ *and with torsion free quotient* $V/U$. *If* $\mu_{L_1}(U) > \mu_{L_1}(V)$, *then* $\mu_{L_0}(U) < \mu_{L_0}(V)$.

*Proof.* We consider the two cases $\mathrm{rank}(U) = 1$ and $\mathrm{rank}(U) = 2$ separately.

Case 1: $\mathrm{rank}(U) = 1$. From $U \hookrightarrow V$, we have either $U \hookrightarrow \mathcal{O}_X(F)$ or $U \hookrightarrow V'$. Note that since $L_1 \cdot \zeta < 0$ and $\zeta = (3F - c_1)$, we have $F \cdot L_1 < \mu_{L_1}(V)$. Now, the first case is impossible since $F \cdot L_1 < \mu_{L_1}(V)$ and $c_1(U) \cdot L_1 > \mu_{L_1}(V)$. Thus, $U \hookrightarrow V'$. Since $V'$ is $L_0$-semistable and no destablizing rank-1 subsheaf of $V'$ can be lifted to $V$, we obtain $\mu_{L_0}(U) < \mu_{L_0}(V') = \mu_{L_0}(V)$.

Case 2: $\mathrm{rank}(U) = 2$. Let $U_2$ be the image of $U$ in $V'$, and let $U_1$ be the kernel of the natural surjection $U \to U_2 \to 0$. Then, $U$ is an extension of $U_2$ by $U_1$. We claim that $U_1 = 0$: otherwise, $U_1$ is a subline bundle of $\mathcal{O}_X(F)$; thus, $(F - c_1(U_1))$ is effective or trivial, and $c_1(U_1) \cdot L_1 \le F \cdot L_1$; since $V'$ is $L_1$-semistable and since $U_2$ is a subsheaf of $V'$, we have $\mu_{L_1}(U_2) \le \mu_{L_1}(V') = (c_1 - F) \cdot L_1/2$; it follows that

$$\mu_{L_1}(U) = \frac{1}{2}[c_1(U_1) \cdot L_1 + c_1(U_2) \cdot L_1] \le \frac{1}{2}[F \cdot L_1 + \frac{1}{2}(c_1 - F) \cdot L_1]$$
$$= \frac{1}{4}(F \cdot L_1 + c_1 \cdot L_1) < \frac{1}{3}(c_1 \cdot L_1) = \mu_{L_1}(V)$$

since $L_1 \cdot \zeta < 0$ and $\zeta = (3F - c_1)$; but this contradicts to our assumption that $\mu_{L_1}(U) > \mu_{L_1}(V)$. Thus, $U_1 = 0$ and $U_2 = U$. It follows that $c_1(V') - c_1(U)$ is effective or trivial. We claim that $c_1(V') - c_1(U)$ is strictly effective: otherwise, $c_1(V') - c_1(U) = 0$; since both $V'$ and $U$ are locally free, $V' \cong U$; by our assumption, $V'$ is $L_0$-semistable; if $V'$ is strictly $L_0$-semistable, then any $L_0$-destablizing rank-1 subsheaf of $V'$ can be lifted to $V$ through $V' \cong U$, contradicting to the assumption; if $V'$ is $L_0$-stable, then the isomorphism between $V'$ and $U$ must be the identity map up to a scalar, and thus the extension (1.17) splits, again contradicting to the assumption. Therefore, $c_1(V') - c_1(U)$ is strictly effective and $\mu_{L_0}(U) < \mu_{L_0}(V')$. Since $L_0 \cdot \zeta = 0$, we see that $\mu_{L_0}(V') = \mu_{L_0}(V)$. Hence $\mu_{L_0}(U) < \mu_{L_0}(V)$. $\square$

Come back to the proof of the proposition 1.16. Since $V$ is $L_1$-unstable, we let

$$0 \subset V_1 \subset \ldots \subset V_k \subset V$$

be the Harder-Narasimhan filtration of $V$ with respect to $L_1$. Assume that $k = 1$. Since $\mu_{L_1}(V_1) > \mu_{L_1}(V)$, we have $\mu_{L_0}(V_1) < \mu_{L_0}(V)$ by the Claim. Put $\zeta = [r \cdot c_1(V_1) - \operatorname{rank}(V_1) \cdot c_1]$. Then by the same argument as in the proof of Lemma 1.3 (ii), we conclude that $\zeta$ defines a wall of type $(3; c_1, c_2)$ separating $L_0$ and $L_1$, and that $W^\zeta$ contains neither $L_0$ nor $L_1$. But this contradicts to the assumption that any wall of type $(3; c_1, c_2)$ which separates $L_0$ and $L_1$ must contain $L_0$.

Finally, we assume that $k = 2$. Then, $\operatorname{rank}(V_1) = 1$ and $\operatorname{rank}(V_2) = 2$. Since $\mu_{L_1}(V_1) > \mu_{L_1}(V)$ and $\mu_{L_1}(V_2/V_1) > \mu_{L_1}(V/V_1)$, one verifies that $\mu_{L_1}(V_2) > \mu_{L_1}(V)$. Thus, we have $\mu_{L_1}(V_i) > \mu_{L_1}(V)$ for $i = 1$ and 2. By the Claim, $\mu_{L_0}(V_i) < \mu_{L_0}(V)$ for $i = 1$ and 2. Now, by similar arguments as in the proof of Lemma 1.4 (ii), we conclude that there exist an integer $r'$ with $0 < r' < 3$ and a divisor $G$ such that $(3G - r'c_1)$ defines a wall of type $(3; c_1, c_2)$ separating $L_0$ and $L_1$, and that the wall $W^{(3G - r'c_1)}$ contains neither $L_0$ nor $L_1$. Again, we obtain a contradiction. $\square$

## 2. Number of moduli of certain bundles.

In this section, $X$ stands for a rational ruled surface. We restrict our attention to rank-3 stable vector bundles on $X$.

### 2.1. The set-up.

Fix a divisor $c_1$ on $X$ and an integer $c_2$ with $(3c_2 - c_1^2) > 0$. Assume that $\mathcal{C}_1$ and $\mathcal{C}_2$ be two adjacent chambers of type $(3; c_1, c_2)$ which share a common face (see Remark 1.10 (i)). Let $L_1$ and $L_2$ be two ample divisors such that $L_1 \in \mathcal{C}_1$ and $L_2 \in \mathcal{C}_2$. Let $V \in \mathfrak{M}_{L_1}(3; c_1, c_2)$. By Proposition 1.14, $V$ is either $L_2$-stable or $L_2$-unstable. Assume that $V$ is $L_2$-unstable. By Theorem 1.6, there exists an integer $i$ with $0 < i < 3$ and a divisor $F$ such that

   (i) $(3F - ic_1) \cdot L_1 < 0 < (3F - ic_1) \cdot L_2$;
   (ii) $-4 \cdot (3c_2 - c_1^2) \le (3F - ic_1)^2 < 0$.

Assume that the Harder-Narasimhan filtration of $V$ with respect to $L_2$ is:

$$0 \subset V_1 \subset \ldots \subset V_k \subset V$$

Then, by Remark 1.7, one of the following is true:

   (a) $k = 1$, $F = c_1(V_1)$, and $i = \operatorname{rank}(V_1) = 1$;
   (b) $k = 1$, $F = c_1(V_1)$, and $i = \operatorname{rank}(V_1) = 2$;
   (c) $k = 2$, $F = c_1(V_2)$, and $i = \operatorname{rank}(V_2) = 2$;
   (d) $k = 2$, $F = c_1(V_1)$, and $i = \operatorname{rank}(V_1) = 1$.

In the following, we shall estimate the number of moduli of vector bundles $V$ which are $L_1$-stable but $L_2$-unstable bundles. Since we have the above four cases, we shall give an formula for the number of moduli case by case in the next four subsections. The basic techniques are those in [11, 12].

We remark that some arithmetic calculations are quite involved. However, the nature of arguments used in these four subsections is quite similar. So one only needs to read proves in subsection 2.2. In fact, if one understands subsection 2.2, one should be able to understand the rest three subsections 2.3, 2.4 and 2.5 without much difficulty since the statements of the results and the arguments in the proofs are similar.

We also emphasize that eventually we only consider two cases for the first Chern class: $c_1 = 0$ and $c_1 = \sigma + f$. One can work out other cases without much difficulty. We single out the case $c_1 = \sigma + f$ because we will consider rank-3 stable vector bundles over the projective plane $\mathbb{P}^2$ in section 3: let $\rho\colon \mathbb{F}_1 \to \mathbb{P}^2$ be the blow-up of $\mathbb{P}^2$ at one point, and let $L$ be the divisor represented by a line in $\mathbb{P}^2$; then $\mathbb{F}_1$ is a rational ruled surface; moreover, $\rho^*(L) = \sigma + f$; up to the dual and up to the twisting by a suitable line bundle, the first Chern class of a rank-3 bundle on $\mathbb{P}^2$ is either zero or $L$; since stable bundles on $\mathbb{P}^2$ are closely related to those on $\mathbb{F}_1$, it follows that it suffices to study rank-3 stable bundles on $\mathbb{F}_1$ with first Chern classes zero or $\sigma + f$ in order to understand an arbitrary rank-3 stable bundle on $\mathbb{P}^2$.

For the sake of clarity, we fix the following notations throughout the rest of this paper. Again, $X$ stands for a rational ruled surface, $\pi\colon X \to \mathbb{P}^1$ is the ruling, $f$ is a fiber of $\pi$, and $\sigma$ is a section to $\pi$ such that $\sigma^2$ is the least. Let $e = -\sigma^2$, and let $K_X$ be the canonical divisor of $X$. Then, $K_X = -2\sigma - (2+e)f$. For an ample divisor $L$ of the form $(x\sigma + yf)$, we let $r_L = y/x$.

## 2.2. Number of moduli of $V$ satisfying case (a) in 2.1.

In this subsection, we estimate the number of moduli of those $V$'s which satisfy the case (a) in subsection 2.1. In this case, $V_1 = \mathcal{O}_X(F)$. Let $W$ be the quotient sheaf $V/V_1$. Then, by Proposition 1.2, $W$ is torsion free and $L_2$-semistable. We have two exact sequences

$$0 \to \mathcal{O}_X(F) \to V \to W \to 0 \qquad \text{and} \qquad 0 \to W \to W^{**} \to Q \to 0. \qquad (2.1)$$

where $Q$ is a sheaf supported on some 0-cycle in $X$. We shall estimate the number of moduli of these $V$'s through the two exact sequences in (2.1). To start with, we estimate the number of moduli of those $Q$, $W^{**}$ and $W$ in the next two lemmas.

**Lemma 2.2.** (i) $\#(moduli\ of\ Q) - dim\ Aut(Q) \leq 2\ell(Q)$;
(ii) $\#(moduli\ of\ W^{**}) - dim\ Aut(W^{**}) \leq 4c_2(W^{**}) - c_1(W^{**})^2 - 4$.

*Proof.* (i) Note that since $Q$ is a sheaf supported on some 0-cycle, we have

$$\#(\text{moduli of } Q) \leq \dim \operatorname{Ext}^1(Q, Q) = 2 \cdot \dim \operatorname{Hom}(Q, Q).$$

Applying $\operatorname{Hom}(., Q)$ to the second exact sequence in (2.1), we get

$$0 \to \operatorname{Hom}(Q, Q) \to \operatorname{Hom}(W^{**}, Q).$$

Thus, $\dim \operatorname{Hom}(Q, Q) \leq 2\ell(Q)$. It follows that

$$\#(\text{moduli of } Q) - \dim \operatorname{Aut}(Q) \leq 2 \cdot \dim \operatorname{Hom}(Q, Q) - \dim \operatorname{Hom}(Q, Q) \leq 2\ell(Q).$$

(ii) Let $c_i' = c_i(W^{**})$ for $i = 1$ and 2. We divide $W^{**}$ into two parts:

$U_s = \{W^{**} |\ W^{**} \text{ is stable with respect to some ample divisor}\}$
$U_{ss} = \{W^{**} |\ W^{**} \text{ is not contained in } U_s\}.$

By the results in [11], #(moduli of $W^{**} \in U_s$) $\leq 4c_2' - (c_1')^2 - 3$.

Next, we consider $W^{**} \in U_{ss}$. Since $W^{**}$ is $L_2$-semistable, it must be strictly $L_2$-semistable. Thus, we have an exact sequence

$$0 \to \mathcal{O}_X(E) \to W^{**} \to \mathcal{O}_X(c_1' - E) \otimes I_Z \to 0$$

where $(2E - c_1') \cdot L_2 = 0$ and $Z$ is a locally complete intersection 0-cycle. It follows that

$$\#(\text{moduli of } W^{**} \in U_{ss}) \leq 4c_2' - (c_1')^2 - 3 + (2 - t)$$

where $t = (4c_2' - (c_1')^2)/4 \geq 0$. Since $W^{**}$ is not contained in $U_s$, $W^{**}$ is not simple by the Proposition 3.2 in [10]; thus, dim $\text{Aut}(W^{**}) \geq 2$. Therefore,

$$\#(\text{moduli of } W^{**} \in U_{ss} \text{ with } t \geq 1) - \text{dim Aut}(W^{**}) \leq 4c_2' - (c_1')^2 - 4.$$

Finally, we assume that $t = 0$. Since $(2E - c_1') \cdot L_2 = 0$, $(2E - c_1')^2 \leq 0$ with equality if and only if $(2E - c_1') = 0$. Since $-(2E - c_1')^2/4 + \ell(Z) = t = 0$, we have $-(2E - c_1')^2 = \ell(Z) = 0$; so $(2E - c_1') = 0$. Since

$$Ext^1(\mathcal{O}_X(c_1' - E), \mathcal{O}_X(E)) = H^1(X, \mathcal{O}_X(2E - c_1')) = H^1(X, \mathcal{O}_X) = 0,$$

$W^{**} = \mathcal{O}_X(E) \oplus \mathcal{O}_X(c_1' - E)$; thus, there is at most one $W^{**}$ in $U_{ss}$ with $t = 0$.  $\square$

**Lemma 2.3.** #(moduli of $W$) $\leq 4c_2(W) - c_1(W)^2 - 3$.

*Proof.* From the second exact sequence in (2.1), we obtain

$$\#(\text{moduli of } W) \leq \#(\text{moduli of } W^{**}) + \#(\text{moduli of } Q) + \text{dim Hom}(W^{**}, Q)$$

$$-\text{dim Aut}(W^{**}) - \text{dim Aut}(Q) + 1.$$

Note that $c_1(W) = c_1(W^{**})$ and $c_2(W) = c_2(W^{**}) + \ell(Q)$. Now, the conclusion follows immediately from Lemma 2.2 (i) and (ii).  $\square$

In the next lemma, we compute the dimensions of some cohomology groups.

**Lemma 2.4.** (i) $h^0(X, W \otimes \mathcal{O}_X(-F)) = 0 = h^0(X, W^{**} \otimes \mathcal{O}_X(2F - c_1))$;
(ii) $h^0(X, W \otimes \mathcal{O}_X(K_X - F)) = h^2(X, W \otimes \mathcal{O}_X(K_X - F)) = 0$;
(iii) Put $\zeta = (3F - c_1)$. Then, $Ext^1(W, \mathcal{O}_X(F))$ has dimension

$$\frac{3c_2 - c_1^2}{3} - 2 - \frac{\zeta^2}{6} + \frac{K_X \cdot \zeta}{2}.$$

*Proof.* (i) Note that $F \cdot L_2 > c_1 \cdot L_2/3$ by the second inequality in the property (i) in subsection 2.1. Since $W$ is $L_2$-semistable, there is no injection from $\mathcal{O}_X(F)$ to $W$. Thus, $h^0(X, W \otimes \mathcal{O}_X(-F))$ must be zero.

Next, we suppose that $h^0(X, W^{**} \otimes \mathcal{O}_X(2F - c_1)) > 0$. Then, the injection $\mathcal{O}_X(c_1 - 2F) \hookrightarrow W^{**}$ induces $\mathcal{O}_X(c_1 - 2F) \otimes I_Z \hookrightarrow W$ where $Z$ is a 0-cycle. Let $U$ be the preimage of $\mathcal{O}_X(c_1 - 2F) \otimes I_Z$ in $V$. Then, $c_1(U) = (c_1 - F)$, and $c_1(U) \cdot L_1/2 < c_1 \cdot L_1/3$ since $V$ is $L_1$-stable. Simplifying this, we obtain that $F \cdot L_1 > c_1 \cdot L_1/3$. But this contradicts to the first inequality in the property (i) in subsection 2.1. Therefore, $h^0(X, W^{**} \otimes \mathcal{O}_X(2F - c_1)) = 0$.

(ii) Since $(-K_X)$ is effective, by the first equality in (i), we have

$$h^0(X, W \otimes \mathcal{O}_X(K_X - F)) \leq h^0(X, W \otimes \mathcal{O}_X(-F)) = 0.$$

Next, note that $h^2(X, W \otimes \mathcal{O}_X(K_X - F)) = h^2(X, W^{**} \otimes \mathcal{O}_X(K_X - F))$; thus, by the Serre duality and the second equality in (i), we obtain

$$h^2(X, W \otimes \mathcal{O}_X(K_X - F)) = h^0(X, W^{**} \otimes \mathcal{O}_X(2F - c_1)) = 0.$$

(iii) Note that $\mathrm{Ext}^1(W, \mathcal{O}_X(F)) \cong H^1(X, W \otimes \mathcal{O}_X(K_X - F))$. Now, the conclusion follows from (ii) and the Riemann-Roch formula. $\square$

The following is the main result in this subsection. It gives an upper bound to the number of moduli of those rank-3 bundles satisfying case (a) in subsection 2.1.

**Proposition 2.5.** *Fix the divisor $F$. Let $\zeta = (3F - c_1)$, and let*

$$d_\zeta(c_1, c_2) = -\frac{3c_2 - c_1^2}{3} + 2 + \frac{\zeta^2}{6} + \frac{K_X \cdot \zeta}{2}.$$

*Then, $\#(\text{moduli of } V) \leq (6c_2 - 2c_1^2 - 8) + d_\zeta(c_1, c_2)$.*

*Proof.* From the first exact sequence in (2.1), we obtain that

$$\#(\text{moduli of } V) \leq \#(\text{moduli of } W) + \dim \mathrm{Ext}^1(W, \mathcal{O}_X(F)) - 1.$$

Therefore, our conclusion follows from Lemma 2.3 and Lemma 2.4 (iii). $\square$

**Lemma 2.6.** *Let $c_1 = a(\sigma + f)$ where $a = 0$ or $1$. Put $\zeta = (3F - c_1)$.*

(i) *If $a = 0$, then $d_\zeta(c_1, c_2) < 0$;*
(ii) *If $a = 1$, then $d_\zeta(c_1, c_2) \leq 0$. Moreover, $d_\zeta(c_1, c_2) = 0$ if and only if either $\zeta = 2\sigma - (3c_2 - 2)f$, or $e = 0$ and $\zeta = -(3c_2 - 2)\sigma + 2f$.*

*Proof.* Let $\zeta = (x\sigma - yf)$. Since $\zeta \cdot L_1 < 0 < \zeta \cdot L_2$, neither $x$ nor $y$ can be zero; moreover, $x$ and $y$ have the same sign.

(i) By Proposition 2.5, we check that

$$d_\zeta(0, c_2) = -c_2 + (2 - 3x) + \frac{3}{2}(ex + 2y)(1 - x)$$

$$= -c_2 + (2 + 3y) + \frac{3}{2}(-x)(ex + 2y + 2 - e).$$

Using the first equality, we see that $d_\zeta(0, c_2) < 0$ when $x > 0$ and $y > 0$; using the second equality, we see that $d_\zeta(0, c_2) < 0$ when $x < 0$ and $y < 0$.

(ii) Since $\zeta = (3F - c_1)$, $x \equiv 2 \pmod 3$ and $y \equiv 1 \pmod 3$. Note that by our assumption, $(3c_2 - c_1^2) > 0$. By Proposition 2.5, we have

$$d_\zeta(c_1, c_2) = -\frac{3c_2 - c_1^2}{3} + 2 + \frac{-x(ex - 3e + 6) - 2y(x - 3)}{6}.$$

First of all, assume that $x > 0$. Then, $x \geq 2$ and $y > 0$. If $x \geq 3$, then $x \geq 5$; it follows that $d_\zeta(c_1, c_2) < 0$. If $x = 2$, then $d_\zeta(c_1, c_2) = [(e + y) - (3c_2 - c_1^2)]/3$; since

$$-4(3c_2 - c_1^2) \leq \zeta^2 = -4(e + y),$$

we see that $d_\zeta(c_1, c_2) \leq 0$ with equality if and only $x = 2$ and $y = (3c_2 - 2)$, that is, $\zeta = 2\sigma - (3c_2 - 2)f$.

Similarly, when $x < 0$, we can show that $d_\zeta(c_1, c_2) \leq 0$ with equality if and only $e = 0$ and $\zeta = -(3c_2 - 2)\sigma + 2f$.   $\square$

*Remark* 2.7. Let $c_1 = \sigma + f$ and $\zeta = 2\sigma - (3c_2 - 2)f$. Then, $F = \sigma + (1 - c_2)f$, $c_1(W) = c_2 f$ and $c_2(W) = 0$. By the second exact sequence in (2.1), $c_2(W^{**}) + \ell(Q) = c_2(W) = 0$; since $W^{**}$ is $L_2$-semistable, $c_2(W^{**}) \geq 0$; thus, $c_2(W^{**}) = \ell(Q) = 0$, and $W$ is locally free. Since $t = [4c_2(W) - c_1(W)^2]/4 = 0$, $W$ can not be $L_2$-stable by the results in [11]. Thus, $W$ is strictly $L_2$-semistable. As in the proof of Lemma 2.2, we conclude that $W = W^{**} = \mathcal{O}_X(c_2/2 \cdot f) \oplus \mathcal{O}_X(c_2/2 \cdot f)$ (so $c_2$ must be even). Therefore, $V$ sits in the exact sequence:

$$0 \to \mathcal{O}_X(\sigma + (1 - c_2)f) \to V \to \mathcal{O}_X(\frac{c_2}{2} \cdot f) \oplus \mathcal{O}_X(\frac{c_2}{2} \cdot f) \to 0. \qquad (2.8)$$

Conversely, let $V$ be the bundle corresponding to an extension class $e = (e_1, e_2)$ in

$$Ext^1(\mathcal{O}_X(c_2/2 \cdot f)^{\oplus 2}, \mathcal{O}_X(\sigma + (1 - c_2)f)) \cong H^1(X, \mathcal{O}_X(\sigma + (1 - 3c_2/2)f))^{\oplus 2};$$

using Proposition 1.16, we see that if $e_1 \neq 0$ and $e_2 \neq 0$, then $V$ is $L'$-stable where $L'$ is contained in the chamber whose upper wall is $W^\zeta$; modulo isomorphism, all such $V$'s are parameterized by a smooth unirational variety of dimension $(6c_2 - 2c_1^2 - 8)$.

## 2.3. Number of moduli of $V$ satisfying case (b) in 2.1.

In this subsection, we estimate the number of moduli of those $V$'s which satisfy case (b) in subsection 2.1. In fact, by considering the dual $V^*$, we shall substantially reduce this case to case (a). By our assumption about $V$, $V$ sits in an exact sequence

$$0 \to V_1 \to V \to \mathcal{O}_X(c_1 - F) \otimes I_Z \to 0 \qquad (2.9)$$

where $V_1$ is $L_2$-semistable. Dualizing (2.9), we obtain

$$0 \to \mathcal{O}_X(F - c_1) \to V^* \to W \to 0$$

where $W$ is torsion free and $W^* = V_1$. Now, $V^*$ is $L_1$-stable but $L_2$-unstable, and the Harder-Narasimhan filtration of $V^*$ with respect to $L_2$ is: $0 \subset \mathcal{O}_X(F - c_1) \subset V^*$. Note that $c_1(V^*) = -c_1$ and $c_2(V^*) = c_2$. Now, we are back to the case (a) subsection 2.1, and can apply the results in subsection 2.2 to $V^*$. In particular, by Proposition 2.5, we have

$$\#(\text{moduli of } V) = \#(\text{moduli of } V^*) \leq (6c_2 - 2c_1^2 - 8) + d_\zeta(-c_1, c_2) \qquad (2.10)$$

where $\zeta = 3(F - c_1) - (-c_1) = (3F - 2c_1)$. We now prove the following Lemma.

**Lemma 2.11.** *Let* $c_1 = a(\sigma + f)$ *where* $a = 0$ *or* $1$. *Put* $\zeta = (3F - 2c_1)$.

   (i) *If* $a = 0$, *then* $d_\zeta(-c_1, c_2) < 0$;

   (ii) *If* $a = 1$, *then* $d_\zeta(-c_1, c_2) \leq 0$. *Moreover,* $d_\zeta(-c_1, c_2) = 0$ *if and only if either* $\zeta = \sigma - (3c_2 - 5)/2 \cdot f$, *or* $e = 0$ *and* $\zeta = -(3c_2 - 5)/2 \cdot \sigma + f$.

*Proof.* The proof of (i) is the same as in the proof of Lemma 2.6 (i). In the following, we prove (ii). Let $\zeta = 3F - 2(\sigma + f) = x\sigma - yf$. Then, by Proposition 2.5,

$$d_\zeta(-c_1, c_2) = -c_2 + \frac{8 - 3x - e}{3} + \frac{(ex + 2y)(3 - x)}{6}.$$

Thus, $d_\zeta(-c_1, c_2) < 0$ if $x \geq 3$. Assume that $x < 3$. Since $x \equiv 1 \pmod 3$, $x = 1$ or $x \leq -2$. Let $x = 1$. Then, $d_\zeta(-c_1, c_2) = -c_2 + (2y + 5)/3$. Note that

$$c_1(V_1) = F = \sigma + \frac{2 - y}{3} \cdot f \quad \text{and} \quad c_2(V_1) \leq c_2 - F \cdot (c_1 - F) = c_2 - \frac{1 + y}{3}.$$

As in the Lemma 1.10 of [11], since $V_1$ is $L_2$-semistable, we have

$$r_{L_2} \leq e + 2c_2(V_1) - \frac{2 - y}{3} \leq e + 2(c_2 - \frac{1 + y}{3}) - \frac{2 - y}{3};$$

on the other hand, since $\zeta \cdot L_2 > 0$, $(e + y) < r_{L_2}$; thus, $2y < (3c_2 - 2)$. Since $y \equiv 2 \pmod 3$, we see that $d_\zeta(-c_1, c_2) < 0$ unless possibly when $2y = (3c_2 - 5)$.

Finally, we assume that $x \leq -2$. Then, $y \leq -1$. Note that the moduli space $\mathfrak{M}_{L_1}(3; c_1, c_2)$ is nonempty and has dimension $(6c_2 - 2c_1^2 - 8)$. Since

$$d_\zeta(-c_1, c_2) = -\frac{3c_2 - c_1^2}{3} + (2 + y) + \frac{x(3e - 6 - ex - 2y)}{6},$$

we conclude that $d_\zeta(-c_1, c_2) < 0$ unless $e = 0$ and $y = -1$. If $e = 0$ and $y = -1$, as in the preceding paragraph, we can show that $d_\zeta(-c_1, c_2) \leq 0$ and that equality holds if and only if $\zeta = -(3c_2 - 5)/2 \cdot \sigma + f$. $\square$

*Remark 2.12..* Let $c_1 = \sigma + f$ and $\zeta = \sigma - (3c_2 - 5)/2 \cdot f$. Since $\zeta = (3F - 2c_1)$, $c_1(V_1) = F = \sigma + (3 - c_2)/2 \cdot f$. ¿From the proof above, we see that

$$e + \frac{3c_2 - 5}{2} = e + y < r_{L_2} \leq e + 2c_2(V_1) - \frac{2 - y}{3};$$

thus, $c_2(V_1) > (1 + y)/3 = (c_2 - 1)/2$, so $c_2(V_1) \geq (c_2 + 1)/2$. On the other hand, from (2.9), we obtain that $c_2(V_1) + \ell(Z) = (c_2 + 1)/2$. It follows that $c_2(V_1) = (c_2 + 1)/2$ and $Z = \emptyset$. Therefore, $V$ sits in the exact sequence

$$0 \to V_1 \to V \to \mathcal{O}_X(\frac{c_2 - 1}{2} \cdot f) \to 0. \tag{2.13}$$

Note that $V_1$ is $L_2$-semistable and that $e + (3c_2 - 5)/2 < r_{L_2} \leq e + (3c_2 - 1)/2$. As in the proof of the Proposition 3.4 in [11], we verify that $V_1$ sits in the exact sequence:

$$0 \to \mathcal{O}_X(\sigma + (1 - c_2)f) \to V_1 \to \mathcal{O}_X(\frac{c_2 + 1}{2} \cdot f) \to 0. \tag{2.14}$$

Now, $V_1/\mathcal{O}_X(\sigma + (1-c_2)f) = \mathcal{O}_X((c_2+1)/2 \cdot f)$ and $V/V_1 = \mathcal{O}_X((c_2-1)/2 \cdot f)$. Therefore, we see that $V/\mathcal{O}_X(\sigma + (1-c_2)f)$ is given by an extension

$$0 \to \mathcal{O}_X(\frac{c_2+1}{2} \cdot f) \to V/\mathcal{O}_X(\sigma + (1-c_2)f) \to \mathcal{O}_X(\frac{c_2-1}{2} \cdot f) \to 0.$$

Since $\operatorname{Ext}^1(\mathcal{O}_X((c_2-1)/2 \cdot f), \mathcal{O}_X((c_2+1)/2 \cdot f)) = 0$, the extension splits. Thus,

$$V/\mathcal{O}_X(\sigma + (1-c_2)f) = \mathcal{O}_X((c_2-1)/2 \cdot f) \oplus \mathcal{O}_X((c_2+1)/2 \cdot f).$$

We conclude that the bundle $V$ sits in an exact sequence:

$$0 \to \mathcal{O}_X(\sigma + (1-c_2)f) \to V \to \mathcal{O}_X(\frac{c_2-1}{2} \cdot f) \oplus \mathcal{O}_X(\frac{c_2+1}{2} \cdot f) \to 0. \quad (2.15)$$

Conversely, using (2.13) and (2.14), we prove the following result.

**Proposition 2.16.** *Let $c_2$ be odd. Then for generic bundles $V_1$ in (2.14), generic extensions $V$ in (2.13) are $L'$-stable where $L'$ is any ample divisor contained in the chamber whose upper wall is $W^{(\sigma - (3c_2-5)/2 \cdot f)}$; moreover, all such $V$'s are parametrized by a smooth rational variety of dimension $(6c_2 - 2c_1^2 - 8)$.*

*Proof.* Put $\zeta = \sigma - (3c_2 - 5)/2 \cdot f$. By the Theorem B in [11], for a fixed integer $c_2$, all the nonempty moduli space $\mathfrak{M}_L(2; \sigma + (3-c_2)/2 \cdot f, (c_2+1)/2)$ are birational; thus, by the Proposition 3.4 in [11], generic bundles $V_1$ in (2.14) are both $L'$-stable and $L_0$-stable, where $L_0 \in W^\zeta$. Dualizing (2.13) and applying Proposition 1.16, we see that for bundles $V_1$ which are stable with respect to both $L'$ and $L_0$, nontrivial extensions $V$ in (2.13) are $L'$-stable.

Moreover, for a fixed $L_0$-stable bundles $V_1$, we see that

$$\operatorname{Hom}(V_1, \mathcal{O}_X(\frac{c_2-1}{2} \cdot f)) \cong H^0(X, V_1^* \otimes \mathcal{O}_X(\frac{c_2-1}{2} \cdot f)) = 0$$

by (2.14); by (2.13), $\dim \operatorname{Hom}(V_1, V) = 1$; thus, all nontrivial extensions $V$ in (2.14) (modulo isomorphisms) are parametrized by

$$\mathbb{P}(\operatorname{Ext}^1(\mathcal{O}_X(\frac{c_2-1}{2} \cdot f), V_1)).$$

Since all the $L_0$-stable $V_1$'s in (2.13) are parametrized by an open dense subset in

$$\mathbb{P}(\operatorname{Ext}^1(\mathcal{O}_X(\frac{c_2+1}{2} \cdot f), \mathcal{O}_X(\sigma + (1-c_2)f))),$$

all such $V$' are parametrized by a smooth rational variety. $\square$

### 2.4. Number of moduli of $V$ satisfying case (c) in 2.1.

In this subsection, we estimate the number of moduli of those $V$'s which satisfy the case (c) in subsection 2.1. Put $F_1 = c_1(V_1)$ and $F_2 = c_1(V_2) = F$. Then, we have two exact sequences

$$0 \to V_2 \to V \to \mathcal{O}_X(c_1 - F_2) \otimes I_{Z_2} \to 0$$
$$0 \to \mathcal{O}_X(F_1) \to V_2 \to \mathcal{O}_X(F_2 - F_1) \otimes I_{Z_1} \to 0. \quad (2.17)$$

¿From the proof of Lemma 1.4, we see that $(2F_1 - F_2) \cdot L_1 < 0 < (2F_1 - F_2) \cdot L_2$. Thus from the second exact sequence in (2.17), we obtain $c_1(V_2)^2 - 4c_2(V_2) \leq (2F_1 - F_2)^2 < 0$. In particular, for fixed $c_1(V_2)$ and $c_2(V_2)$, there are at most finitely many choices of $F_1$. We start with two lemmas which deals with the dimensions of certain cohomology groups.

**Lemma 2.18.** (i) $h^i(X, \mathcal{O}_X(K_X + F_2 - 2F_1) \otimes I_{Z_1}) = 0$ for $i = 0$ and $2$;

(ii) *The dimension of* $\text{Ext}^1(\mathcal{O}_X(F_2 - F_1) \otimes I_{Z_1}, \mathcal{O}_X(F_1))$ *is equal to*

$$\ell(Z_1) - 1 - (2F_1 - F_2)(2F_1 - F_2 - K_X)/2;$$

(iii) $\#(\text{moduli of } V_2) \leq 3\ell(Z_1) - 2 - (2F_1 - F_2)(2F_1 - F_2 - K_X)/2.$

*Proof.* (i) Note that $\pm(F_2 - 2F_1)$ can not be effective. Thus, we have

$$h^0(X, \mathcal{O}_X(K_X + F_2 - 2F_1) \otimes I_{Z_1}) \leq h^0(X, \mathcal{O}_X(F_2 - 2F_1)) = 0,$$

and $h^2(X, \mathcal{O}_X(K_X + F_2 - 2F_1) \otimes I_{Z_1}) = h^2(X, \mathcal{O}_X(K_X + F_2 - 2F_1)) = 0.$

(ii) By the Serre duality, we obtain the following

$$\text{Ext}^1(\mathcal{O}_X(F_2 - F_1) \otimes I_{Z_1}, \mathcal{O}_X(F_1)) \cong H^1(X, \mathcal{O}_X(K_X + F_2 - 2F_1) \otimes I_{Z_1}).$$

By (i) and the Riemann-Roch formula, our conclusion follows.

(iii) From the second exact sequence in (2.17), we see that

$$\#(\text{moduli of } V_2) \leq \#(\text{moduli of } Z_1) + \dim \text{Ext}^1(\mathcal{O}_X(F_2 - F_1) \otimes I_{Z_1}, \mathcal{O}_X(F_1)) - 1.$$

Therefore, the conclusion follows from (ii).  $\square$

**Lemma 2.19.** (i) $h^0(X, V_2 \otimes \mathcal{O}_X(K_X + c_1 - 2F_2) \otimes I_{Z_2}) = 0;$

(ii) $h^2(X, V_2 \otimes \mathcal{O}_X(K_X + c_1 - 2F_2) \otimes I_{Z_2}) = 0;$

(iii) *Let* $\eta = 3F_2 - 2c_1$. *Then,* $\dim \text{Ext}^1(\mathcal{O}_X(c_1 - F_2) \otimes I_{Z_2}, V_2)$ *is equal to*

$$2c_2 - \frac{5c_1^2}{9} - c_2(V_2) - 2 + \frac{K_X \cdot \eta}{2} - \frac{\eta^2}{18} + \frac{c_1 \cdot \eta}{9}.$$

*Proof.* (i) Since $(-K_X)$ is effective, it suffices to show that

$$h^0(X, V_2 \otimes \mathcal{O}_X(c_1 - 2F_2)) = 0.$$

Suppose the contrary. Tensoring the second exact sequence in (2.17) by $\mathcal{O}_X(c_1 - 2F_2)$, we see that either $(c_1 - 2F_2 + F_1)$ or $(c_1 - F_2 - F_1)$ is effective. On the other hand, we have $\mu_{L_2}(V_2/\mathcal{O}_X(F_1)) > \mu_{L_2}(V/\mathcal{O}_X(F_1))$ and $\mu_{L_2}(\mathcal{O}_X(F_1)) > \mu_{L_2}(V)$ since the Harder-Narasimhan filtration of $V$ with respect to $L_2$ is:

$$0 \subset \mathcal{O}_X(F_1) \subset V_2 \subset V.$$

It follows that $(c_1 - 2F_2 + F_1) \cdot L_2 < 0$ and $(c_1 - F_2 - F_1) \cdot L_2 < 0$. A contradiction.

(ii) We need only to show that $h^0(X, V_2 \otimes \mathcal{O}_X(F_2 - c_1)) = 0$. Suppose the contrary. Then, we have $\mathcal{O}_X(c_1 - F_2) \hookrightarrow V_2 \hookrightarrow V$. Since $V$ is $L_1$-stable, $(c_1 - F_2) \cdot L_1 < c_1 \cdot L_1/3$, so $F_2 \cdot L_1/2 > c_1 \cdot L_1/3$; but this is absurd since $V_2 \hookrightarrow V$.

(iii) Note that by the Serre duality, we have

$$\text{Ext}^1(\mathcal{O}_X(c_1 - F_2) \otimes I_{Z_2}, V_2) \cong H^1(X, V_2 \otimes \mathcal{O}_X(K_X + c_1 - 2F_2) \otimes I_{Z_2}).$$

Therefore, our conclusion follows from (i), (ii) and the Riemann-Roch formula.  $\square$

**Proposition 2.20.** *Fix $F = F_2 = c_1(V_2)$. Let $\eta = (3F - 2c_1)$, and let*

$$d_\eta^*(c_1, c_2) = -\frac{2(3c_2 - c_1^2)}{3} + 3 + \frac{(2F_1 - F_2)^2}{4} + \frac{(2F_1 - F_2) \cdot K_X}{2} + \frac{\eta^2}{12} + \frac{\eta \cdot K_X}{2}.$$

*Then, $\#(moduli\ of\ V) \leq (6c_2 - 2c_1^2 - 8) + d_\eta^*(c_1, c_2)$.*

*Proof.* By the first exact sequence in (2.17), $\#(moduli\ of\ V)$ is bounded above by

$$\#(\text{moduli of } V_2) + \#(\text{moduli of } Z_2) + \dim \text{Ext}^1(\mathcal{O}_X(c_1 - F_2) \otimes I_{Z_2}, V_2) - 1.$$

Therefore, the conclusion follows from Lemma 2.18 (iii) and Lemma 2.19 (iii).  □

Next, we estimate $d_\eta^*(c_1, c_2)$. The argument is similar to but much complicated than those in the proofs of Lemma 2.6 and Lemma 2.11.

**Lemma 2.21.** *Let $c_1 = a(\sigma + f)$ where $a = 0$ or $1$. Put $\eta = (3F - 2c_1)$.*

(i) *If $a = 0$, then $d_\eta^*(c_1, c_2) < 0$;*

(ii) *If $a = 1$, then $d_\eta^*(c_1, c_2) < 0$ unless possibly either of the following:*

    (iia) $\eta = \sigma - (3c_2 - 2)/2 \cdot f$ *or* $\eta = \sigma - (3c_2 - 5)/2 \cdot f$

    (iib) $e = 0$, $\eta = -(3c_2 - 2)/2 \cdot \sigma + f$ *or* $\eta = -(3c_2 - 5)/2 \cdot \sigma + f$.

*Proof.* (i) Put $(2F_1 - F_2) = x_1\sigma - y_1 f$ and $F_2 = x_2\sigma - y_2 f$. Since

$$(x_i\sigma - y_i f) \cdot L_1 < 0 < (x_i\sigma - y_i f) \cdot L_2,$$

neither $x_i$ nor $y_i$ can be zero; moreover, $x_i$ and $y_i$ have the same sign. Note that $c_2 \geq 1$. By Proposition 2.20, $d_\eta^*(0, c_2) = -2c_2 + 3 + d_1 + 3d_2$ where

$$d_1 = \frac{(2F_1 - F_2)^2}{4} + \frac{(2F_1 - F_2) \cdot K_X}{2} = \frac{-x_1(ex_1 - 2e + 4) - 2y_1(x_1 - 2)}{4}$$

$$d_2 = \frac{F_2^2}{4} + \frac{F_2 \cdot K_X}{2} = \frac{-x_2(ex_2 - 2e + 4) - 2y_2(x_2 - 2)}{4}.$$

First of all, we assume that $e \geq 1$. Then, $d_i \leq -2$ if $x_i \neq 1$. When $x_1 = 1$, $d_1 = (e + 2y_1)/4 - 1$; by the first exact sequence in (2.17), $c_2(V_2) \leq (c_2 + F_2^2)$; by the second exact sequence in (2.17), $c_2(V_2) \geq F_1(F_2 - F_1)$; thus,

$$(4F_1^2 - 4F_1 \cdot F_2) \geq -4c_2(V_2) \geq -4(c_2 + F_2^2),$$

and $-(e + 2y_1) = (2F_1 - F_2)^2 \geq -4c_2 - 3F_2^2$; it follows that

$$d_1 \leq (c_2 + 3F_2^2)/4 - 1 < (c_2 - 1).$$

Similarly, when $x_2 = 1$, $d_2 = (e + 2y_2)/4 - 1$; since $-(e + 2y_2) = F_2^2 \geq -4c_2/3$, we have $3d_2 \leq (c_2 - 3)$. In any case, we see that $(-c_2 + d_1) \leq -2$ and $(-c_2 + 3d_2) \leq -2$. Thus, $d_\eta^*(0, c_2) < 0$.

Next, we assume that $e = 0$. Then, $d_i = (-x_i y_i - 2x_i + 2y_i)/2$. Thus, $d_i \leq -2$ if $x_i \geq 2$ or $y_i \leq -2$. When $x_1 = -1$ or $y_1 = -1$, as in the preceding paragraph, we can show that $d_1 < (c_2 - 1)$. Similarly, when $x_2 = -1$ or $y_2 = -1$, we can show that $3d_2 \le (c_2 - 3)$. Thus, $d^*_\eta(0, c_2) < 0$ in all the cases.

(ii) Put $(2F_1 - F_2) = x_1\sigma - y_1 f$ and $\eta = x_2\sigma - y_2 f$. Again, neither $x_i$ nor $y_i$ can be zero; moreover, $x_i$ and $y_i$ have the same sign. Since $\eta = (3F - 2c_1)$, $x_2 \equiv 1 \pmod 3$ and $y_2 \equiv 2 \pmod 3$. Note that $(3c_2 - c_1^2) > 0$. By Proposition 2.20, $d^*_\eta(c_1, c_2) = -2(3c_2 - c_1^2)/3 + 3 + d_1 + d'_2$ where $d_1$ is as in the proof of (i) and

$$d'_2 = \frac{\eta^2}{12} + \frac{\eta \cdot K_X}{2} = \frac{-x_2(ex_2 - 6e + 12) - 2y_2(x_2 - 6)}{12}.$$

As in the proof of (i), $d_1 - (3c_2 - c_1^2)/3 \le \eta^2/12 - 1$. First of all, we assme that $x_2 > 1$. Then, either $x_2 \ge 7$ or $x_2 = 4$. If $x_2 \ge 7$, then $d'_2 \le -7$, and $d^*_\eta(c_1, c_2) < 0$. If $x_2 = 4$, then $d'_2 = (2e + y_2)/3 - 4$. Since $-4(3c_2 - c_1^2) \le \eta^2 = -8(2e + y_2)$, we see that $(2e + y_2) \le (3c_2 - c_1^2)/2$ and that $d'_2 \le (3c_2 - c_1^2)/6 - 4$; again, $d^*_\eta(c_1, c_2) < 0$.

Next, we assume that $x_2 = 1$. By the assumption made in subsection 2.1, $W^\eta$ is the only wall of type $(3; c_1, c_2)$ separating $L_1$ and $L_2$. Let $L_0 \in W^\eta$. Since $V$ is $L_1$-stable, $V$ is $L_0$-semistable by Theorem 1.6 and Remark 1.10 (ii); as in the proof of the Theorem 3.1 in [9], we must have $r_{L_0} < e + (3c_2 + 1)/2$; since $L_0 \cdot \eta = 0$, $r_{L_0} = e + y_2$, so $y_2 < (3c_2 + 1)/2$. Since $2y_2 \equiv 1 \pmod 3$, $2y_2 \le (3c_2 - 2)$. Note that $d'_2 = 5(e + 2y_2)/12 - 1$. If $x_1 \ne 1$, then $d_1 \le -2$, and $d^*_\eta(c_1, c_2) < 0$. If $x_1 = 1$, then $d_1 - (3c_2 - c_1^2)/3 \le \eta^2/12 - 1 = -(e + 2y_2)/12 - 1$; thus,

$$d^*_\eta(c_1, c_2) \le -\frac{3c_2 - c_1^2}{3} + 3 + (-\frac{e + 2y_2}{12} - 1) + (\frac{5(e + 2y_2)}{12} - 1) = \frac{2y_2 + 5 - 3c_2}{3}.$$

It follows that $d^*_\eta(c_1, c_2) < 0$ unless possibly when $y_2 = (3c_2 - 2)/2$ or $(3c_2 - 5)/2$.

Finally, we assume that $x_2 < 0$. When $e \ge 1$, we can show that $d'_2 \le -2$; thus, $d^*_\eta(c_1, c_2) < 0$. When $e = 0$, as in the preceding paragraph, $d^*_\eta(c_1, c_2) < 0$ unless possibly when $y_2 = -1$ and $x_2 = -(3c_2 - 2)/2$ or $-(3c_2 - 5)/2$. $\square$

Remark 2.22.  $c_1 = \sigma + f$ and $\eta = \sigma - (3c_2 - 2)/2 \cdot f$. Since $\eta = 3F - 2c_1$,

$$c_1(V_2) = F = \sigma + (2 - c_2)/2 \cdot f.$$

By the first exact sequence in (2.17), $c_2(V_2) + \ell(Z_2) = c_2/2$. On the other hand, let $L_0 \in W^\eta$; then, $V$ is $L_0$-semistable; since $c_1(V_2) \cdot L_0/2 = c_1 \cdot L_0/3$, $V_2$ must be $L_0$-semistable; thus, by the Lemma 1.10 of [11], we obtain

$$e + \frac{3c_2 - 2}{2} = r_{L_0} \le e + 2c_2(V_2) + \frac{c_2 - 2}{2};$$

so $c_2(V_2) \ge c_2/2$, and it follows that $c_2(V_2) = c_2/2$ and $Z_2$ is empty. As in Remark 2.12, we verify that $V_2$ sits in an exact sequence

$$0 \to \mathcal{O}_X(\sigma + (1 - c_2)f) \to V_2 \to \mathcal{O}_X(\frac{c_2}{2} \cdot f) \to 0,$$

and that $V$ sits in the exact sequence

$$0 \to \mathcal{O}_X(\sigma + (1 - c_2)f) \to V \to \mathcal{O}_X(\frac{c_2}{2} \cdot f) \oplus \mathcal{O}_X(\frac{c_2}{2} \cdot f) \to 0.$$

We conclude that the bundles $V$ here are exactly those in subsection 2.2.

*Remark* 2.23. $c_1 = \sigma + f$ and $\eta = \sigma - (3c_2 - 5)/2 \cdot f$. From the proof of (ii), $d_\eta^*(c_1, c_2) \leq 0$. Assume that $d_\eta^*(c_1, c_2) = 0$. Then,

$$d_1 - (3c_2 - c_1^2)/3 = \eta^2/12 - 1 = -(e + 3c_2 - 5)/12 - 1;$$

since $d_1 = (e + 2y_1)/4 - 1$, we see that $y_1 = (3c_2 - 1)/2$. Calculating the second Chern classes of $V_2$ and $V$ by (2.17), we conclude that both $Z_1$ and $Z_2$ must be empty. Now, the two exact sequences in (2.17) are simplified into:

$$0 \to V_2 \to V \to \mathcal{O}_X(\frac{c_2 - 1}{2} \cdot f) \to 0$$

$$0 \to \mathcal{O}_X(\sigma + (1 - c_2)f) \to V_2 \to \mathcal{O}_X(\frac{c_2 + 1}{2} \cdot f) \to 0.$$

These are exactly the same as (2.13) and (2.14). Therefore, the bundles $V$ here are exactly those bundles $V$ in subsection 2.3.

## 2.5. Number of moduli of $V$ satisfying case (d) in 2.1.

In this subsection, we estimate the number of moduli of those $V$'s which satisfy the last case (d) in subsection 2.1. Put $F_1 = c_1(V_1) = F$ and $F_2 = c_1(V_2)$. Then, we have

$$0 \to \mathcal{O}_X(F_1) \to V \to V/\mathcal{O}_X(F_1) \to 0$$

$$0 \to V_2/\mathcal{O}_X(F_1) \to V/\mathcal{O}_X(F_1) \to \mathcal{O}_X(c_1 - F_2) \otimes I_Z \to 0. \quad (2.25)$$

¿From the proof of Lemma 1.4, we see that

$$(2F_2 - F_1 - c_1) \cdot L_1 < 0 < (2F_2 - F_1 - c_1) \cdot L_2.$$

Dualizing the exact sequences in (2.25), we obtain

$$0 \to [V/\mathcal{O}_X(F_1)]^* \to V^* \to \mathcal{O}_X(-F_1) \otimes I_{Z_2} \to 0$$

$$0 \to \mathcal{O}_X(F_2 - c_1) \to [V/\mathcal{O}_X(F_1)]^* \to \mathcal{O}_X(F_1 - F_2) \otimes I_{Z_1} \to 0. \quad (2.26)$$

Now, $V^*$ is $L_1$-stable but $L_2$-unstable. Note that

$$2(F_2 - c_1) - c_1([V/\mathcal{O}_X(F_1)]^*) = (2F_2 - F_1 - c_1).$$

Thus, as in the preceding subsection, we conclude that

$$\#(\text{moduli of } V) = \#(\text{moduli of } V^*) \leq (6c_2 - 2c_1^2 - 8) + d_\eta^*(-c_1, c_2) \quad (2.27)$$

where $\eta = 3c_1([V/\mathcal{O}_X(F_1)]^*) - 2c_1(V^*) = (3F_1 - c_1)$, and $d_\eta^*(-c_1, c_2)$ is equal to

$$-\frac{2(3c_2 - c_1^2)}{3} + 3 + \frac{(2F_2 - F_1 - c_1)^2}{4} + \frac{(2F_2 - F_1 - c_1) \cdot K_X}{2} + \frac{\eta^2}{12} + \frac{\eta \cdot K_X}{2}. \quad (2.28)$$

We omit the proof of the following lemma which is an analogue of Lemma 2.21.

**Lemma 2.29.** *Let $c_1 = a(\sigma + f)$ where $a = 0$ or $1$. Put $\eta = (3F_1 - c_1)$.*

(i) *If $a = 0$, then $d_\eta^*(-c_1, c_2) < 0$;*

(ii) *If $a = 1$, then $d_\eta^*(-c_1, c_2) < 0$ unless possibly when either $\eta = 2\sigma - (3c_2 - 2)f$, or $e = 0$ and $\eta = -(3c_2 - 2)f + 2f$.*

*Remark* 2.30. Let $c_1 = \sigma + f$ and $\eta = 2\sigma - (3c_2 - 2)f$. Since $\eta = (3F_1 - c_1)$, we obtain $F_1 = \sigma + (1 - c_2)f$; thus, $c_1(V/\mathcal{O}_X(F_1)) = c_2 f$ and $c_2(V/\mathcal{O}_X(F_1)) = 0$. Let $L_0 \in W^\eta$. Then, $V$ is $L_0$-semistable and $c_1 \cdot L_0/3 = c_1(V/\mathcal{O}_X(F_1)) \cdot L_0/2$. Thus, $V/\mathcal{O}_X(F_1)$ is $L_0$-semistable. As in Remark 2.7, we see that $V/\mathcal{O}_X(F_1)$ is equal to $\mathcal{O}_X(c_2/2 \cdot f)^{\oplus 2}$. Therefore, $V$ sits in the exact sequence:

$$0 \to \mathcal{O}_X(\sigma + (1 - c_2)f) \to V \to \mathcal{O}_X(\frac{c_2}{2} \cdot f) \oplus \mathcal{O}_X(\frac{c_2}{2} \cdot f) \to 0.$$

We conclude that the bundles $V$ here are exactly those in subsection 2.2.

## 3. Moduli space sof rank-3 stable bundles.

In this section, $X$ stands for a rational ruled surface. We will combine the results in [9] and in section 2 of this paper to study the irreducibility, unirationality and rationality of the moduli space of rank-3 bundles stable with respect to an arbitrary ample divisor. As an application, we study the moduli space of rank-3 stable bundles on the projective plane $\mathbb{P}^2$.

### 3.1. The moduli space $\mathfrak{M}_L(3; 0, c_2)$.

**Theorem 3.1.** *Let $L_1$ and $L_2$ be two ample divisors on a rational ruled surface. Then,*

(i) *The moduli space $\mathfrak{M}_{L_1}(3; 0, c_2)$ is nonempty if and only if $c_2 \geq 3$;*

(ii) *If $c_2 \geq 3$, then $\mathfrak{M}_{L_1}(3; 0, c_2)$ and $\mathfrak{M}_{L_2}(3; 0, c_2)$ are birational. In particular, they are smooth, irreducible and unirational with dimension $(6c_2 - 8)$.*

*Proof.* Note that since $(-K_X)$ is effective, $\mathfrak{M}_{L_i}(3; 0, c_2)$ is smooth with dimension $(6c_2 - 8)$ whenever it is nonempty. We may assume that $r_{L_1} \gg 0$ such that $L_1$ does not lie on any wall of type $(3; 0, c_2)$. By the Theorem 5.4 and Remark 5.7 in [9], $\mathfrak{M}_{L_1}(3; 0, c_2)$ is nonempty if and only if $c_2 \geq 3$; moreover, if $c_2 \geq 3$, then $\mathfrak{M}_{L_1}(3; 0, c_2)$ is smooth, irreducible and unirational with dimension $(6c_2 - 8)$.

Next, we assume that $L_2$ does not lie on any wall of type $(3; 0, c_2)$. When $L_2$ lies in the same chamber as $L_1$ does, then by Theorem 1.6, $\mathfrak{M}_{L_2}(3; 0, c_2)$ is naturally identified with $\mathfrak{M}_{L_1}(3; 0, c_2)$; thus, the conclusion follows. Assume that $L_2$ and $L_1$ lie in different chambers. Since the set of walls of type $(3; 0, c_2)$ is locally finite, we can pick up finitely many $\mathbb{Q}$-ample divisors

$$L_1 = L^{(1)}, L^{(2)}, \ldots, L^{(k-1)}, L^{(k)} = L_2$$

on the line segment connecting $L_1$ and $L_2$ such that every $L^{(i)}$ lies in some chamber, and that $L^{(i)}$ and $L^{(i+1)}$ are separated by a single wall for each $1 \leq i \leq (k - 1)$. For each $1 \leq i \leq (k - 1)$, we apply the set-up in subsection 2.1 to $L^{(i)}$ and $L^{(i+1)}$; by Lemma 2.6 (i), Lemma 2.11 (i), Lemma 2.21 (i) and Lemma 2.29 (i), $\mathfrak{M}_{L^{(i)}}(3;0,c_2)$ and $\mathfrak{M}_{L^{(i+1)}}(3;0,c_2)$ are identified by removing some subschemes of codimension at least one; it follows that $\mathfrak{M}_{L_1}(3;0,c_2)$ and $\mathfrak{M}_{L_2}(3;0,c_2)$ are identified by removing some subschemes of codimension at least one. It follows that $\mathfrak{M}_{L_2}(3;0,c_2)$ is smooth, irreducible and unirational.

Finally, we assume that $L_2$ lies on some wall $W$ of type $(3;0,c_2)$. We claim that $\mathfrak{M}_{L_2}(3;0,c_2)$ is nonempty: let $H_1$ and $H_2$ be two ample divisors not lying on any wall of type $(3;0,c_2)$ such that $L_2$ lies on the line segment connecting $H_1$ and $H_2$; from the preceding paragraph, there exists $V \in \mathfrak{M}_{H_1}(3;0,c_2)$ such that $V$ is also $H_2$-stable; since $L_2$ is a linearly combination of $H_1$ and $H_2$ with positive coefficients, $V$ is also $L_2$-stable. Thus, $\mathfrak{M}_{L_2}(3;0,c_2)$ is nonempty, and has dimension $(6c_2 - 8)$. Let $\mathcal{C}$ be a chamber whose closure contains $L_2$, and let $L'_2 \in \mathcal{C}$. By Corollary 1.15, $\mathfrak{M}_{L_2}(3;0,c_2)$ is naturally embedded in $\mathfrak{M}_{L'_2}(3;0,c_2)$. By the preceding paragraph, $\mathfrak{M}_{L'_2}(3;0,c_2)$ is smooth, irreducible and unirational with dimension $(6c_2 - 8)$. Since $\mathfrak{M}_{L_2}(3;0,c_2)$ and $\mathfrak{M}_{L'_2}(3;0,c_2)$ have the same dimension and since $\mathfrak{M}_{L'_2}(3;0,c_2)$ is irreducible, it follows that $\mathfrak{M}_{L_2}(3;0,c_2)$ is smooth, irreducible and unirational. $\qquad\square$

### 3.2. The moduli space $\mathfrak{M}_L(3;\sigma + f, c_2)$.

**Theorem 3.2.** *Assume that the integer $c_2$ is even.*

(i) *If $c_2 < 2$, then all the moduli spaces $\mathfrak{M}_L(3;\sigma + f, c_2)$ are empty;*

(ii) *Assume that $c_2 \geq 2$. If $r_L \geq e + (3c_2 - 2)/2$ or $r_L \leq 2/(3c_2 - 2)$, then $\mathfrak{M}_L(3;\sigma + f, c_2)$ is empty. If $2/(3c_2 - 2) < r_L < e + (3c_2 - 2)/2$, then $\mathfrak{M}_L(3;\sigma + f, c_2)$ is irreducible, smooth, and unirational with dimension $(6c_2 + 2e - 12)$; moreover, a generic bundle in $\mathfrak{M}_L(3;\sigma + f, c_2)$ sits in an exact sequence*

$$0 \to \mathcal{O}_X(\sigma + (1 - c_2)f) \to V \to \mathcal{O}_X(\tfrac{c_2}{2} \cdot f) \oplus \mathcal{O}_X(\tfrac{c_2}{2} \cdot f) \to 0. \qquad (3.3)$$

*Proof.* (i) By the Theorem 1.1 in [9], $\mathfrak{M}_H(3;\sigma + f, c_2)$ is empty if $r_H \gg 0$. If $\mathfrak{M}_L(3;\sigma + f, c_2)$ is nonempty for some $L$, then all bundles in $\mathfrak{M}_L(3;\sigma + f, c_2)$ are not $H$-stable; thus, generic bundles in an irreducible component of $\mathfrak{M}_L(3;\sigma + f, c_2)$ must come from some nonempty wall $W^\zeta$ as in subsection 2.1, and $d_\zeta(\sigma + f, c_2) \geq 0$ or $d^*_\zeta(\sigma + f, c_2) \geq 0$. Since $c_2$ is even, by Lemma 2.6 (ii), Lemma 2.11 (ii), Lemma 2.21 (ii) and Lemma 2.29 (ii), we have either of the following two cases:

(a) $\zeta = 2\sigma - (3c_2 - 2)f$ or $\sigma - (3c_2 - 2)/2 \cdot f$.

(b) $e = 0$, $\zeta = -(3c_2 - 2)\sigma + 2f$ or $-(3c_2 - 2)/2 \cdot \sigma + f$.

Since $c_2 \leq 0$, the class $\zeta$ is effective in all the above cases; but this is absurd since we must have $\zeta \cdot L_0 < 0$ for some ample divisor $L_0$.

(ii) Case (iia): $e \geq 1$. Since $r_L > e \geq 1$ and $c_2 \geq 2$, we always have $r_L > 2/(3c_2-2)$. If $r_L \geq e+(3c_2-2)/2$, then $L \cdot (2\sigma - (3c_2-2)f) \geq 0$; as in the proof of (i), we verify that $\mathfrak{M}_L(3;\sigma + f, c_2)$ is empty. Next, we assume that $r_L < e + (3c_2-2)/2$. Let $L'$ be contained in the chamber whose upper wall is $W^{(2\sigma-(3c_2-2)f)}$. By Lemma 2.6 (ii), Lemma 2.11 (ii), Lemma 2.21 (ii) and Lemma 2.29 (ii), if $W^\eta$ separates $L'$ and $L$, then the number of moduli of bundles coming from $W^\eta$ is less than $(6c_2+2e-12)$; thus, $\mathfrak{M}_{L'}(3;\sigma + f, c_2)$ and $\mathfrak{M}_L(3;\sigma + f, c_2)$ are identified by removing some subschemes of codimension at least one. By Remark 2.7, $\mathfrak{M}_{L'}(3;\sigma + f, c_2)$ is irreducible and unirational, and a generic bundle in $\mathfrak{M}_{L'}(3; \sigma + f, c_2)$ sits in (3.3). Therefore, the conclusion for the moduli space $\mathfrak{M}_L(3; \sigma + f, c_2)$ follows.

Case (iib): $e = 0$. Note that there is a symmetry between $\sigma$ and $f$. By Case (iia), $\mathfrak{M}_L(3; \sigma + f, c_2)$ is empty if $r_L \geq (3c_2 - 2)/2$ or $r_L \leq 2/(3c_2 - 2)$. If

$$2/(3c_2 - 2) < r_L < (3c_2 - 2)/2,$$

then again, we can show the conclusion as we did in Case (iia).   $\square$

**Theorem 3.4.** *Assume that the integer $c_2$ is odd.*

(i) *If $c_2 < 3$, then all the moduli spaces $\mathfrak{M}_L(3; \sigma + f, c_2)$ are empty;*

(ii) *Assume that $c_2 \geq 3$. If $r_L \geq e + (3c_2 - 5)/2$ or $r_L \leq 2/(3c_2 - 5)$, then $\mathfrak{M}_L(3; \sigma + f, c_2)$ is empty. If $2/(3c_2 - 5) < r_L < e + (3c_2 - 5)/2$, then $\mathfrak{M}_L(3; \sigma + f, c_2)$ is irreducible, smooth, and rational with dimension $(6c_2 + 2e - 12)$; moreover, a generic bundle in $\mathfrak{M}_L(3; \sigma + f, c_2)$ sits in an exact sequence*

$$0 \to \mathcal{O}_X(\sigma + (1 - c_2)f) \to V \to \mathcal{O}_X(\frac{c_2 - 1}{2} \cdot f) \oplus \mathcal{O}_X(\frac{c_2 + 1}{2} \cdot f) \to 0. \quad (3.5)$$

*Proof.* The proof of (i) is similar to that of Theorem 3.2 (i). For (ii), let $L'$ be contained in the chamber whose upper wall is $W^{(\sigma - (3c_2 - 5)/2 \cdot f)}$; by Proposition 2.16 and (2.15), $\mathfrak{M}_{L'}(3; \sigma + f, c_2)$ is irreducible, smooth, and rational with dimension $(6c_2 + 2e - 12)$, and a generic bundle in $\mathfrak{M}_{L'}(3; \sigma + f, c_2)$ sits in (3.5); now, the rest of the proof is similar to the proof of Theorem 3.2 (ii).   $\square$

### 3.3. Moduli of stable rank-3 bundles on $\mathbb{P}^2$.

Let $\rho : \mathbb{F}_1 \to \mathbb{P}^2$ be the blow-up morphism of $\mathbb{P}^2$ at one point, and let $\sigma$ be the exceptional divisor. Then, $\mathbb{F}_1$ is the rational ruled surface with $e = 1$ and $\sigma$ is the section to the ruling such that $\sigma^2 = -e = -1$. Let $L$ be the divisor in $\mathbb{P}^2$ represented by a line. Then, $\rho^* L = (\sigma + f)$ where $f$ is a fiber of the ruling on $\mathbb{F}_1$. Put $L_n = n\rho^* L - \sigma$. If $n \gg 0$, then $L_n$ is an ample divisor on $\mathbb{F}_1$. As in [6, 2], we have the following (we omit its proof since the proof is identical to that in [6]).

**Lemma 3.6.** *For $n \gg 0$, there is an open immersion*

$$\phi : \mathfrak{M}_L(3; c_1, c_2) \hookrightarrow \mathfrak{M}_{L_n}(3; \rho^* c_1, c_2) \quad (3.7)$$

*which is induced by $\phi(V) = \rho^*(V)$ on closed points.*

In [3], Drezet and Le Potier have determined the nonemptyness of the moduli space $\mathfrak{M}_L(3; c_1, c_2)$, and showed that $\mathfrak{M}_L(3; c_1, c_2)$ is irreducible whenever it is nonempty. For instance, by the Theorem B in [3], we verify that $\mathfrak{M}_L(3; 0, c_2)$ is nonempty if and only if $c_2 \geq 3$, and that $\mathfrak{M}_L(3; L, c_2)$ is nonempty if and only if $c_2 \geq 2$. The following improves upon Drezet and Le Potier's result when $r = 3$.

**Theorem 3.8.** *If $c_2 \geq 3$, then $\mathfrak{M}_L(3; 0, c_2)$ is irreducible, smooth, and unirational with dimension $(6c_2 - 8)$. If $c_2 \geq 2$, then $\mathfrak{M}_L(3; L, c_2)$ is irreducible, smooth, and unirational (rational when $c_2$ is odd) with dimension $(6c_2 - 10)$.*

*Proof.* Note that if $\mathfrak{M}_L(3; c_1, c_2)$ is nonempty, then it is smooth with dimension $(6c_2 - 2c_1^2 - 8)$ by a result of Maruyama. Using the open immersion (3.7), we see that the results follow from Theorem 3.1, Theorem 3.2 and Theorem 3.4.   $\square$

DEPARTMENT OF MATHEMATICS, HKUST, CLEAR WATER BAY, KOWLOON, HONG KONG
 *E-mail address*: mawpli@masu1.ust.hk

DEPARTMENT OF MATHEMATICS, OKLAHOMA STATE UNIVERSITY, STILLWATER, OK 74078, USA
 *E-mail address*: zq@math.okstate.edu